\documentclass[12pt,a4paper,english,reqno]{amsart}
\usepackage[a4paper,footskip=1.5em]{geometry}
\usepackage{amsmath,amssymb,amsthm,mathtools,bbm}
\usepackage[mathscr]{euscript}
\usepackage[usenames,dvipsnames]{color}
\usepackage{adjustbox,tikz,calc,graphics,babel,standalone}
\usepackage{subcaption}
\usepackage{csquotes,enumerate,verbatim}
\usepackage[final]{microtype}
\usepackage[numbers]{natbib}
\usetikzlibrary{shapes.misc,calc,intersections,patterns,decorations.pathreplacing}
\usepackage{hyperref}
\hypersetup{colorlinks=true,linkcolor=blue,citecolor=blue,pdfpagemode=UseNone,pdfstartview={XYZ null null 1.00}}
\usepackage{cmtiup}

\allowdisplaybreaks

\theoremstyle{plain}
\newtheorem*{theorem*}{Theorem}
\newtheorem{theorem}{Theorem}[section]
\newtheorem{lemma}[theorem]{Lemma}
\newtheorem{claim}[theorem]{Claim}
\newtheorem{proposition}[theorem]{Proposition}
\newtheorem*{claim*}{Claim}

\theoremstyle{remark}
\newtheorem*{remark}{Remark}

\def\N{\mathbb{N}}

\def\R{\mathbb{R}}
\def\P{\mathbb{P}}

\def\C{\mathcal}

\def\Scr{\mathscr}
\def\Cond{\,|\,}
\DeclareMathOperator\Bi{Bin}
\DeclareMathOperator\Po{Po}

\let\eps\varepsilon

\let\originalleft\left
\let\originalright\right
\renewcommand{\left}{\mathopen{}\mathclose\bgroup\originalleft}
\renewcommand{\right}{\aftergroup\egroup\originalright}

\def\HE{{\widehat E}}
\def\AM{{\mathbf{M}}}
\def\omegap{{\omega'}}

\makeatletter
\def\imod#1{\allowbreak\mkern10mu({\operator@font mod}\,\,#1)}
\makeatother

\begin{document}

\title{Line percolation}

\author{Paul Balister}
\address{Department of Mathematical Sciences, University of Memphis, Memphis TN 38152, USA}
\email{pbalistr@memphis.edu}

\author{B\'{e}la Bollob\'{a}s}
\address{Department of Pure Mathematics and Mathematical Statistics, University of Cambridge, Wilberforce Road, Cambridge CB3\thinspace0WB, UK, {\em and\/}
Department of Mathematical Sciences, University of Memphis, Memphis TN 38152, USA, {\em and\/} London Institute for Mathematical Sciences, 35a South St., Mayfair, London W1K\thinspace2XF, UK}
\email{b.bollobas@dpmms.cam.ac.uk}

\author{Jonathan Lee}
\address{Department of Pure Mathematics and Mathematical Statistics, University of Cambridge, Wilberforce Road, Cambridge CB3\thinspace0WB, UK}
\email{j.d.lee@dpmms.cam.ac.uk}

\author{Bhargav Narayanan}
\address{Department of Pure Mathematics and Mathematical Statistics, University of Cambridge, Wilberforce Road, Cambridge CB3\thinspace0WB, UK}
\email{b.p.narayanan@dpmms.cam.ac.uk}

\date{21 February 2015}

\begin{abstract}
We study a new geometric bootstrap percolation model, \emph{line percolation}, on the $d$-dimensional integer grid $[n]^d$. In line percolation with infection parameter $r$, infection spreads from a subset $A\subset [n]^d$ of initially infected lattice points as follows: if there exists an axis-parallel line $L$ with $r$ or more infected lattice points on it, then every lattice point of $[n]^d$ on $L$ gets infected, and we repeat this until the infection can no longer spread. The elements of the set $A$ are usually chosen independently, with some density $p$, and the main question is to determine $p_c(n,r,d)$, the density at which percolation (infection of the entire grid) becomes likely. In this paper, we determine $p_c(n,r,2)$ up to a multiplicative factor of $1+o(1)$ and $p_c(n,r,3)$ up to a multiplicative constant as $n\rightarrow \infty$ for every fixed $r\in \N$. We also determine the size of the minimal percolating sets in all dimensions and for all values of the infection parameter.
\end{abstract}

\maketitle

\section{Introduction}
Bootstrap percolation models and arguments have been used to study a range of phenomena in various areas, ranging from crack formation and the dynamics of glasses to neural nets and economics; see~\citep{i3,i1,i2} for a small sample of such applications. In this paper, we introduce and study a new geometric bootstrap percolation model defined on the $d$-dimensional integer grid $[n]^d$ which we call \emph{line percolation}; here, we write $[n]$ for the set $\{1, 2, \dots, n$\}. For $v \in [n]^d$, let $\C{L}(v)$ denote the set of $d$ axis-parallel lines through $v$ and let
\[\C{L}(n,d) = \bigcup_{v\in [n]^d}\C{L}\left(v\right)\]
denote the set of all axis-parallel lines that pass through the points of $[n]^d$. In line percolation with infection parameter $r\in\N$ (or $r$-neighbour line percolation for short), infection spreads from a set of initially infected lattice points as follows: if there exists a line with $r$ or more infected lattice points on it, then every lattice point on that line becomes infected. More precisely, given a set $A \subset [n]^d$ of \emph{initially infected points}, we have a sequence $A = A^{(0)} \subset A^{(1)} \subset \dots \subset A^{(t)} \subset \dots$ of subsets of $[n]^d$
such that
\[A^{\left(t+1\right)} = A^{\left(t\right)} \cup \left\{ v \in [n]^d : \exists L \in \C{L}\left(v\right) \mbox{ such that } |L \cap A^{\left(t\right)}| \ge r\right\}.\]
We say that a point $v \in [n]^d$ is \emph{infected} at time $t$ if $v \in A^{(t)}$, and we say that a line $L \in \C{L}(n,d)$ is \emph{active} at time $t$ if $L \subset A^{(t)}$. The \emph{closure} of $A$ is the set $[A] = \bigcup_{t \ge 0} A^{(t)}$ of eventually infected points. We say that the process \emph{terminates} when no more newly infected points are added. If every point of $[n]^d$ is infected (or equivalently, if every line is active) when the process terminates, i.e., if $[A] = [n]^d$, then we say that $A$ \emph{percolates}.

The classical model of \emph{$r$-neighbour bootstrap percolation on a graph} was introduced by Chalupa, Leath and Reich~\citep{Chalupa79} in the context of disordered magnetic systems and has since been extensively studied not only by mathematicians but also by physicists and sociologists; see~\citep{Adler03, ising1, socio1, socio2}, for example. In this model, a vertex of the graph gets infected if it has at least $r$ previously infected neighbours in the graph. The model is usually studied in the random setting where the main question is to determine the critical threshold at which percolation occurs. If the elements of the initially infected set are chosen independently at random, each with probability $p$, then one aims to determine the value $p_c$ at which percolation becomes likely. In this regard, the $r$-neighbour bootstrap percolation model on $[n]^d$, with edges induced by the integer lattice, has been the subject of a large body of work; see~\citep{Holroyd03, Balogh09, Balogh12}, and the references therein.

On account of its inherent geometric structure, it is possible to construct other interesting bootstrap percolation models on the $d$-dimensional grid. In the past, this has involved endowing the grid with a graph structure other than the one induced by the integer lattice (which, in other words, is a Cartesian product of paths).

In this direction, Holroyd, Liggett and Romik~\citep{cross_nbd} considered $r$-neighbour bootstrap percolation on $[n]^2$ where the neighbourhood of a lattice point $v$ is taken to be a `cross' centred at $v$ consisting of $r-1$ points in each of the four axis directions. Sharp thresholds for a model with an anisotropic variant of these cross neighbourhoods were obtained recently by Duminil-Copin and van Enter~\citep{anisotropic}. Gravner, Hoffman, Pfeiffer and Sivakoff~\citep{hamming} studied the $r$-neighbour bootstrap percolation model on the $d$-dimensional Hamming torus which is the graph on $[n]^d$ where $u,v \in [n]^d$ are adjacent if and only if $u-v$ has exactly one nonzero coordinate; the Hamming torus, in other words, is the Cartesian product of complete graphs, which is perhaps the second most natural graph structure on $[n]^d$ after the grid. They obtained bounds on the critical exponents (i.e., $\log_n (p_c)$) which are tight for all values of the infection parameter when $d=2$ and for a few small values of the infection parameter when $d=3$.

The line percolation model we consider here is a natural variant of the bootstrap percolation model on the Hamming torus studied by Gravner, Hoffman, Pfeiffer and Sivakoff. However, we should note that while all the other models mentioned above are $r$-neighbour bootstrap percolation models on \emph{some underlying graph}, the line percolation model \emph{is not}.

Morally, line percolation is better thought of as coming from the very general \emph{$\C{U}$-bootstrap percolation} model introduced by Bollob\'as, Smith and Uzzell~\citep{nf1}. In this model, one starts by specifying a finite collection $\C{U}$ of finite subsets of an infinite lattice; a point $v$ of the lattice becomes infected when, for some $U \in \C{U}$, all the points in the translate $U+v$ of the set $U$ are infected. In their paper, Bollob\'as, Smith and Uzzell prove a classification theorem for two-dimensional models of this type and show that every such model is either \emph{supercritical}, \emph{critical} or \emph{subcritical}. In particular, they show that a model is supercritical if and only if there exist finite subsets of the lattice from which the infection can spread to the whole lattice. While line percolation on the integer lattice cannot be described by associating a \emph{finite} family of neighbourhoods with each point of the lattice, there do exist, as we shall see, finite sets from which the infection can spread to the whole lattice, and our results about the critical probabilities of the line percolation model are in agreement with the general bounds for the critical probabilities of supercritical models proved in~\citep{nf1}. For some related work concerning subcritical models, see the paper of Balister, Bollob\'as, Przykucki and Smith~\citep{nf2}.

\section{Our results}
In this note, our main aim is to investigate what happens in the line percolation model when the initial set $A_p\subset [n]^d$ of infected points is determined by randomly selecting points from $[n]^d$, each independently with probability $p \in (0,1)$. We shall primarily be concerned with the following question: for what values of $p$ is percolation likely to occur? Let $\vartheta_p(n,r,d)$ denote the probability that such a randomly chosen initial set $A_p$ percolates. We note that $\vartheta_p(n,r,d)$ increases with $p$ and define the \emph{critical probability} $p_c(n,r,d)$ by setting
\[p_c\left(n,r,d\right) = \inf \left\{p :  \vartheta_p\left(n,r,d\right) \ge 1/2 \right\}.\]
The primary question of interest is to determine the asymptotic behaviour of $p_c(n,r,d)$ as $n \rightarrow \infty$ for fixed  $d,r \in \N$. Note that when $r=1$, a set of lattice points percolates if and only if it is nonempty, so $p_c(n,1,d) = \Theta(n^{-d})$; therefore, we restrict our attention to the case where $r \ge 2$.

Before we state our results, a few remarks about asymptotic notation are in order. We shall make use of standard asymptotic notation; the variable tending to infinity will always be $n$ unless we explicitly specify otherwise. When convenient, we shall also make use of some notation (of Vinogradov) that might be considered non-standard: given functions $f(n)$ and $g(n)$, we write $f \ll g$ if $f = O(g)$, $f \gg g$ if $g = O(f)$, and $f \sim g$ if $f=(1+o(1))g$. Constants suppressed by the asymptotic notation are allowed to depend on the fixed parameters ($c \in (0,1)$ and $d,r \in \N$), but not on $n$ or $p$.

In two dimensions, we are able to estimate the probability of percolation $\vartheta_p(n,r,2)$ up to a multiplicative constant for all $0\le p \le 1$. We also determine $p_c(n,r,2)$ up to a multiplicative factor of $1+o(1)$ as $n \rightarrow \infty$.

\begin{theorem}\label{2d-unboosted-perc}
Fix $r, s \in \N$ with $r \ge 2$ and $0 \le s \le r-1$. Then as $n \to \infty$,
\begin{equation} \vartheta_p\left(n,r,2\right) = \Theta\left(n^{2s+1}\left(np\right)^{r\left(2s+1\right) - s\left(s + 1\right)}\right) \label{2dformula}
\end{equation}
when $n^{-1-\frac{1}{r-s-1}} \ll p \ll n^{-1-\frac{1}{r-s}}$, with the convention that $n^{-1-\frac{1}{r-s-1}} = 0$ when $s=r-1$. Also, $\vartheta_p(n,r,2) = \Theta(1)$ when $p \gg n^{-1-\frac{1}{r}}$. Furthermore, we have
\[ p_c\left(n,r,2\right) \sim \lambda n^{-1-\frac{1}{r}},\]
where $\lambda = (r!\log 2/2)^{1/r}$.
\end{theorem}

The techniques used to obtain the above formula for $\vartheta_p(n,r,2)$ allow us to prove the following result about the critical probability in three dimensions, which is the main result of this paper.

\begin{theorem}\label{3d-perc}
Fix an integer $r \ge 2$, and let $s = \lfloor \sqrt{r+ 1/4} - 1/2\rfloor$. As $n \to \infty$, writing $\gamma = (r + s^2+s)/(2s+2)$, we have
\[ p_c\left(n,r,3\right) = \Theta\left(n^{-1-\frac{1}{r-\gamma}}\right).\]
\end{theorem}

The nature of the threshold at the critical probability is also worth investigating. We say that the model exhibits a sharp threshold at $p=p(n,r,d)$ if for any fixed $\eps > 0$, we have $\vartheta_{(1+\eps)p}(n,r,d) = 1 - o(1)$ and $\vartheta_{(1-\eps)p}(n,r,d) = o(1)$. It is not difficult to see from Theorem~\ref{2d-unboosted-perc} and our proof of Theorem~\ref{3d-perc} that in stark contrast to the classical $r$-neighbour bootstrap percolation model on the grid, there is \emph{no sharp threshold phenomenon} in the line percolation model in dimensions two and three. We expect similar behaviour in higher dimensions but we do not have a proof of such an assertion.

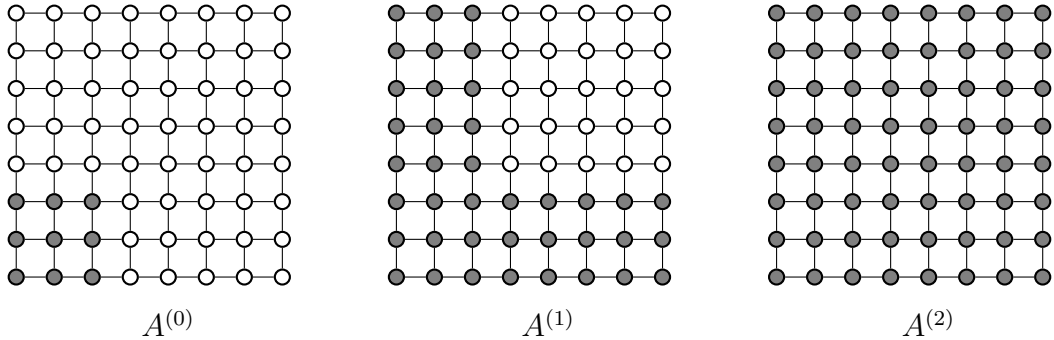
\begin{figure}
\begin{center}
\begin{tikzpicture}
\foreach \x in {4,5,6,7,8}
\foreach \y in {4,5,6,7,8}
	\node (l\x\y) at (\x/2-5, \y/2) [inner sep=0.7mm, thick, circle, draw=black!100, fill=black!0] {};

\foreach \x in {1,2,3}
\foreach \y in {4,5,6,7,8}
	\node (l\x\y) at (\x/2-5, \y/2) [inner sep=0.7mm, thick, circle, draw=black!100, fill=black!0] {};

\foreach \y in {1,2,3}
\foreach \x in {4,5,6,7,8}
	\node (l\x\y) at (\x/2-5, \y/2) [inner sep=0.7mm, thick, circle, draw=black!100, fill=black!0] {};

\foreach \x in {1,2,3}
\foreach \y in {1,2,3}
	\node (l\x\y) at (\x/2-5, \y/2) [inner sep=0.7mm, thick, circle, draw=black!100, fill=black!50] {};	

\foreach \x in {4,5,6,7,8}
\foreach \y in {4,5,6,7,8}
	\node (m\x\y) at (\x/2, \y/2) [inner sep=0.7mm, thick, circle, draw=black!100, fill=black!0] {};

\foreach \x in {1,2,3}
\foreach \y in {4,5,6,7,8}
	\node (m\x\y) at (\x/2, \y/2) [inner sep=0.7mm, thick, circle, draw=black!100, fill=black!50] {};

\foreach \y in {1,2,3}
\foreach \x in {4,5,6,7,8}
	\node (m\x\y) at (\x/2, \y/2) [inner sep=0.7mm, thick, circle, draw=black!100, fill=black!50] {};

\foreach \x in {1,2,3}
\foreach \y in {1,2,3}
	\node (m\x\y) at (\x/2, \y/2) [inner sep=0.7mm, thick, circle, draw=black!100, fill=black!50] {};

\foreach \x in {1,2,3,4,5,6,7,8}
\foreach \y in {1,2,3,4,5,6,7,8}
	\node (r\x\y) at (\x/2+5, \y/2) [inner sep=0.7mm, thick, circle, draw=black!100, fill=black!50] {};	

\foreach \x in {1,2,3,4,5,6,7,8}
\foreach \y[evaluate={\yp=int(\y+1)}] in {1,2,3,4,5,6,7}
	\draw (l\x\y) -- (l\x\yp);
	
\foreach \x[evaluate={\xp=int(\x+1)}] in {1,2,3,4,5,6,7}
\foreach \y in {1,2,3,4,5,6,7,8}
	\draw (l\x\y) -- (l\xp\y);

\foreach \x in {1,2,3,4,5,6,7,8}
\foreach \y[evaluate={\yp=int(\y+1)}] in {1,2,3,4,5,6,7}
	\draw (m\x\y) -- (m\x\yp);
	
\foreach \x[evaluate={\xp=int(\x+1)}] in {1,2,3,4,5,6,7}
\foreach \y in {1,2,3,4,5,6,7,8}
	\draw (m\x\y) -- (m\xp\y);
	
\foreach \x in {1,2,3,4,5,6,7,8}
\foreach \y[evaluate={\yp=int(\y+1)}] in {1,2,3,4,5,6,7}
	\draw (r\x\y) -- (r\x\yp);
	
\foreach \x[evaluate={\xp=int(\x+1)}] in {1,2,3,4,5,6,7}
\foreach \y in {1,2,3,4,5,6,7,8}
	\draw (r\x\y) -- (r\xp\y);	

\node at (5/2 - 5,-0.1) {$A^{(0)}$};	
\node at (5/2,-0.1) {$A^{(1)}$};
\node at (5/2 + 5,-0.1) {$A^{(2)}$};	
\end{tikzpicture}
\end{center}
\caption{The spread of infection from $A=[3]^2$ in the $3$-neighbour line percolation process on $[8]^2$.}
\label{minsetfig}
\end{figure}

It is also an interesting question to determine the size of a minimal percolating set for $r$-neighbour line percolation on $[n]^d$ for any $d, r \in \N$ and $n \ge r+1$. It is easy to check that the set $[r]^d$ percolates (see Figure~\ref{minsetfig}). We shall demonstrate that this is in fact optimal.

\begin{theorem}\label{min-set}
The minimum size of a percolating set in the $r$-neighbour line percolation model on $[n]^d$ is $r^d$ for all $d,r, n\in\N$ with $n \ge r+1$.
\end{theorem}

Establishing this fact turns out to be harder than it appears at first glance. The result is trivial when $d = 1$. When $d = 2$, it is not hard to demonstrate that any percolating set has size at least $r^2$. Indeed, consider a generalised two-dimensional line percolation model on $[n]^2$ where the infection thresholds for horizontal and vertical lines are $r_h$ and $r_v$ respectively; note that we recover the $r$-neighbour line percolation model when $r_h = r_v = r$. Let $M(r_h, r_v)$ denote the size of a minimal percolating set in this generalised model. If the first line $L$ to be infected is horizontal, then $L$ must contain $r_h$ initially infected points and furthermore, if the set of initially infected points is a percolating set, then the set of initially infected points not on $L$ must constitute a percolating set for the generalised process with infection parameters $r_h$ and $r_v - 1$. An analogous statement holds if $L$ is vertical, so it follows that
\[M(r_h, r_v) \ge \min \left( r_v + M(r_h - 1, r_v), r_h + M(r_h, r_v-1)\right).\]
We obtain by induction that $M(r_h, r_v) \ge r_h r_v$ which implies in particular that $M(r, r) \ge r^2$. The argument described above depends crucially on the fact that a line has codimension one in a two-dimensional space. The incidence geometry of a collection of lines in the plane is essentially straightforward; this is no longer the case in higher dimensions and we need more delicate arguments to prove Theorem~\ref{min-set}.

The rest of this paper is organised as follows. We collect together some useful facts in Section~\ref{binvar}. We consider line percolation in two dimensions in Section~\ref{2d} and prove Theorem~\ref{2d-unboosted-perc}. In Section~\ref{3d}, we turn to line percolation in three dimensions and prove Theorem~\ref{3d-perc}. We give the proof of Theorem~\ref{min-set} in  Section~\ref{minset}. We conclude the paper in Section~\ref{remark} with some discussion.

\section{Preliminaries}\label{binvar}
We will need some standard facts about binomial random variables. We collect these here for the sake of convenience. As is usual, for a random variable with distribution $\Bi(N,p)$, we write $\mu = Np$ for its mean.

The first proposition we shall require is an easy consequence of the fact that $e^{-2x} \le 1 - x \le e^{-x}$ for all $0 \le x \le 1/2$.

\begin{proposition}
\label{binsmall}
Let $X$ be a random variable with distribution $\Bi (N,p)$ and suppose that $p\le 1/2$. For every $1 \le k\le N$, we have
\[ \exp \left(-2 \mu \right) \left(\mu/k\right)^k \le \P \left(X = k\right) \le \exp \left(- \mu \right) \left(2e\mu/k\right)^k.\]
Also, $ \exp (-2 \mu ) \le \P (X = 0) \le \exp (- \mu )$. \qed
\end{proposition}

We shall make use of  the following standard concentration result; see~\citep{probtextbook} for example.

\begin{proposition}
\label{chernoff}
Let $X$ be a random variable with distribution $\Bi(N,p)$. For any $0<\delta<1$, we have
\[ \P\left( |X - \mu| \ge \delta \mu \right) \le 2 \exp \left(\frac{-\delta^2 \mu}{3} \right). \eqno\qed \]
\end{proposition}

Finally, we shall make use of the following easy proposition.

\begin{proposition}
\label{smallmu}
Let $X$ be a random variable with distribution $\Bi(N,p)$. For any fixed $k\ge 0$, as $N \to \infty$,
\[ \P\left( X \ge k\right) = \begin{cases}
       \Theta \left(\P\left( X = k\right)\right) = \Theta (\mu^k) \hfill & \text{ if } \mu \ll 1, \text{ and}\\
       \Theta(1) \hfill & \text{ if } \mu \gg 1.
  \end{cases}
\]
In particular, for any fixed $k\ge 0$, we always have $\P( X \ge k) = O(\mu^k)$. \qed
\end{proposition}

The FKG inequality of Fortuin, Kasteleyn and Ginibre~\citep{FKG} asserts that a pair of events are positively correlated under a product measure if they are both increasing, and that they are negatively correlated if one of these events is increasing and the other is decreasing. Recall that an event $E \subset \{0,1\}^{[n]}$ is \emph{increasing} if $\omega \in E$ implies that $\omegap \in E$ for every $\omegap \in \{0,1\}^{[n]}$ such that $\omegap_x \ge \omega_x$ for each $x\in [n]$; \emph{decreasing} events are defined analogously. We shall need a simple corollary of the FKG inequality. Given a set $A \subset [n]$, we say that an event $E \subset \{0,1\}^{[n]}$ is \emph{decreasing on $A$} if $\omega \in E$ implies that $\omegap \in E$ for every $\omegap \in \{0,1\}^{[n]}$ such that $\omegap_x \le \omega_x$ for each $x\in A$ and $\omegap_x = \omega_x$ for each $x \notin A$.

\begin{proposition}
\label{fkg}
Let $A \subset [n]$ and let $\P$ be a product measure on $\{0,1\}^{[n]}$. For any increasing event $F$ which depends only on the coordinates in $A$ and any event $E$ which is decreasing on $A$, we have $\P(F\,|\,E) \le \P(F)$.
\end{proposition}
\begin{proof}
For $v \in \{0,1\}^{[n]\setminus A}$, denote by $I_v$ the event that the coordinates in $[n] \setminus A$ are given by $v$. Since $E$ is decreasing on $A$ and $F$ is increasing on $A$, we see by applying the FKG inequality to the induced product measure on $\{0,1\}^A$ that
\[\P(E\cap F \Cond I_v) \le \P(E\Cond I_v)\P(F\Cond I_v)\]
for every $v \in \{0,1\}^{[n]\setminus A}$. Since $F$ does not depend on the coordinates in $[n] \setminus A$, we also have $\P(F\Cond I_v) = \P(F)$ for every $v \in \{0,1\}^{[n]\setminus A}$. Therefore, by summing over all $v \in \{0,1\}^{[n]\setminus A}$, we see that
\begin{align*}
\P(E \cap F) &= \sum_{v} \P (I_v)\P(E \cap F \Cond  I_v)\\
&\le \sum_{v} \P (I_v)\P(E\Cond I_v)\P(F\Cond I_v)\\
&= \sum_{v}\P (I_v) \P(E\Cond I_v)\P(F) = \P(E)\P(F). \qedhere
\end{align*}
\end{proof}

\section{Line percolation in two dimensions}\label{2d}
The proof of the following proposition is essentially identical to the proof of Theorem 2.1 in~\citep{hamming}; we reproduce it here for completeness.

\begin{proposition}\label{2d-magnitude}
Fix $r\in \N$ with $r \ge 2$, and let $C > 0$ be a positive constant. If $p = C n^{-1-{1/r}}$, then
\[ \vartheta_p(n,r,2) \sim 1 - \exp{(-2C^r/r!)}.\]
\end{proposition}
\begin{proof}
Let $E_h$ denote the event that there exists a horizontal line containing $r$ or more initially infected points and define $E_v$ analogously. Clearly, the process terminates on the first step if neither $E_h$ nor $E_v$ holds; so $\vartheta_p \le \P (E_h \cup E_v)$.

The number of horizontal lines containing at least $r$  initially infected points is binomially distributed and it is easily seen to converge in distribution to a $\Po(C^r / r!)$ random variable. Therefore, $\P(E_h) \sim 1 - \exp(-C^r/r!)$ and similarly, $\P(E_v) \sim 1 - \exp(-C^r/r!)$.

We now estimate $\P (E_h \cap E_v)$. Let $E_h \circ E_v$ denote the event that $E_h$ and $E_v$ occur disjointly. Now, $E_h$ and $E_v$ are increasing events, so it follows from the FKG inequality~\citep{FKG} and the BKR inequality~\citep{BKR-1, BKR-2} that $\P (E_h \cap E_v) \ge \P(E_h)\P(E_v) \ge \P (E_h \circ E_v)$. Observe that $(E_h \cap E_v) \setminus (E_h \circ E_v)$ occurs only if some lattice point $v$ is initially infected and each of the two axis-parallel lines through $v$ contains at least $r-1$ initially infected points. It follows that
\[ \P \left((E_h \cap E_v) \setminus (E_h \circ E_v)\right) = O\left(n^2 p (np)^{2r-2} \right) = O\left(n^{-1+1/r}\right), \]
so $\P ((E_h \cap E_v) \setminus (E_h \circ E_v)) = o(1)$. Consequently, we see that $\P(E_h \cap E_v) \sim \P(E_h)\P(E_v)$. Hence,
\[ \P(E_h \cup E_v) \sim \P(E_h) + \P(E_v) - \P(E_h)\P(E_v) \sim 1 - \exp{(-2C^r/r!)}, \]
so $\vartheta_p \le 1 - \exp{(-2C^r/r!)} + o(1)$.

To bound $\vartheta_p$ from below, we generate the initial configuration in two rounds, first by sprinkling infected points with density $p' = p (1-1/\log n)$ and then (independently) with density $p'' = p/\log n$; clearly, this configuration is dominated by an initial configuration where points are infected independently with density $p$, so it suffices to bound from below the probability that such an initial configuration percolates.

Let $E'$ be the event that there exists a line that contains at least $r$ initially infected points from the \emph{first sprinkling} of points. It is easy to check from the argument above that $\P (E') \sim 1 - \exp(-2C^r/r!)$.

Now, fix a line $L$ and note that the probability that a particular line perpendicular to $L$ has at least $r-1$ initially infected points (none of which are on $L$) from the \emph{second sprinkling} of points is $\Theta( ((n-1)p'')^{r-1} ) = \Theta(n^{-1+1/r}(\log n)^{-r+1})$. Thus, the number of such lines perpendicular to $L$ is a binomial random variable with mean $\mu = \Omega(n^{1/r}(\log n)^{-r})$. Since $\mu \rightarrow \infty$ as $n \rightarrow \infty$, by Proposition~\ref{chernoff}, the probability that there exist at least $r$ such lines perpendicular to $L$ in the second sprinkling is $1-o(1)$. Hence, conditional on $E'$, the probability of percolating using the points infected in the second sprinkling is $1-o(1)$. Therefore,
\[ \vartheta_p \ge (1-o(1))\P(E')  = 1 - \exp{(-2C^r/r!)} - o(1)\]
and the result follows.
\end{proof}

We shall now prove Theorem~\ref{2d-unboosted-perc}.
\begin{proof}[Proof of Theorem~\ref{2d-unboosted-perc}]
It follows from Proposition~\ref{2d-magnitude} that $p_c(n,r,2) \sim \lambda n^{-1-{1/r}}$, where $\lambda$ is the unique positive real number such that $\exp (-2\lambda^r/ r!) = 1/2$. This implies that if $p \gg n^{-1-1/r}$, then 	$\vartheta_p(n,r,2) = \Theta(1)$.

We now deal with the case where $p \ll n^{-1-1/r}$. Fix an integer $s \in \{ 0,1,\dots,r-1\}$ and suppose that
\[n^{-1-\frac{1}{r-s-1}} \ll p \ll n^{-1-\frac{1}{r-s}},\] with the convention that $n^{-1-1/(r-s-1)}$ is equal to $0$ when $s = r-1$. Observe that $n(np)^{r-i} \ll 1$ for each $0 \le i \le s$ and that $n(np)^{r-i} \gg 1$ for each $ i \ge s+1$.
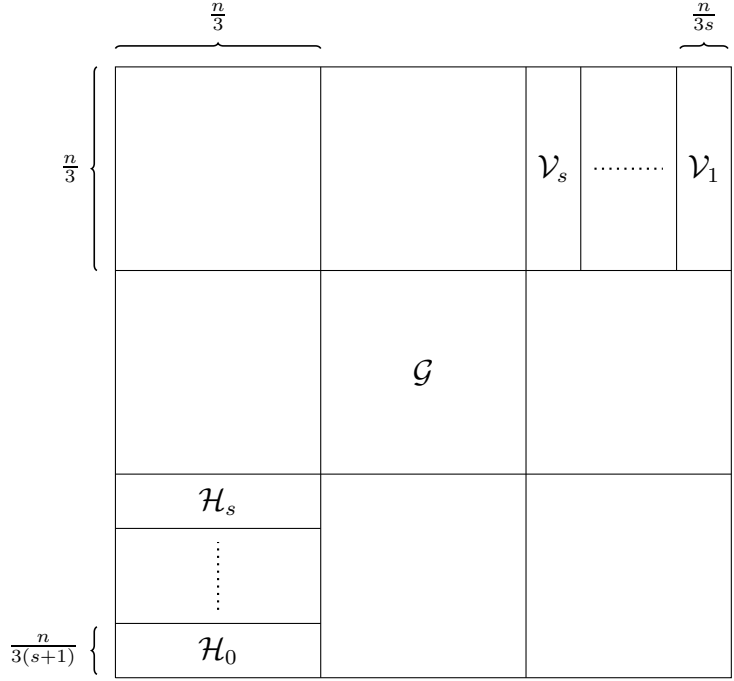
\begin{figure}
\begin{center}
\begin{tikzpicture}[scale=0.9]
\clip (-1.8,-0.5) rectangle (11,10);
\draw (0,0)--(9,0)--(9,9)--(0,9)--(0,0);
\draw (0,3)--(9,3);
\draw (0,6)--(9,6);
\draw (3,0)--(3,9);
\draw (6,0)--(6,9);
\draw [decorate,decoration={brace,amplitude=2pt}, semithick, yshift=2pt]
(0,9.2) -- (3,9.2) node [black,midway,yshift=12pt]
{\footnotesize $\tfrac{n}{3}$};
\draw [decorate,decoration={brace,amplitude=2pt}, semithick, xshift=-2pt]
(-0.2, 6) -- (-0.2, 9) node [black,midway,xshift=-10pt]
{\footnotesize $\tfrac{n}{3}$};
\node (s) at (4.5,4.5) [] {$\C{G}$};

\draw (0,0.8)--(3,0.8);
\draw (0,2.2)--(3,2.2);
\draw [dotted, thick] (1.5, 1.0) -- (1.5,2.0);
\node (h1) at (1.5,0.4) [] {$\C{H}_0$};
\node (h2) at (1.5,2.6) [] {$\C{H}_{s}$};

\draw [decorate,decoration={brace,amplitude=2pt}, semithick, xshift=-2pt]
(-0.2, 0.05) -- (-0.2, 0.75) node [black,midway,xshift=-20pt]
{\footnotesize $\tfrac{n}{3(s+1)}$};

\draw (6.8,6)--(6.8,9);
\draw (8.2,6)--(8.2,9);
\draw [dotted, thick] (7, 7.5) -- (8,7.5);
\node (v1) at (8.6,7.5) [] {$\C{V}_1$};
\node (h2) at (6.4,7.5) [] {$\C{V}_{s}$};

\draw [decorate,decoration={brace,amplitude=2pt}, semithick, yshift=2pt]
(8.25,9.2) -- (8.95,9.2) node [black,midway,yshift=12pt]
{\footnotesize $\tfrac{n}{3s}$};

\end{tikzpicture}
\end{center}
\caption{The lower bound construction.}
\label{constr}
\end{figure}

We first bound $\vartheta_p(n,r,2)$ from below. Let $\C{H}_0, \C{H}_1, \dots, \C{H}_{s}$, $\C{G}$ and $\C{V}_1, \C{V}_2, \dots, \C{V}_s$ be disjoint subgrids of $[n]^2$ as in Figure~\ref{constr}. We bound the probability of percolation from below in terms of the following events. For $0 \le i \le s$, let $E_h(i)$ denote the event that there exists a horizontal line $L$ passing through the $n/(3s+3) \times n/3$ `horizontal' grid $\C{H}_i$ such that there exist at least $r-i$ initially infected points in $ L \cap \C{H}_i$. Analogously, for $1 \le i \le s$, let $E_v(i)$ denote the event that there exists a vertical line $L$ passing through the $n/3 \times n/3s$ `vertical' grid $\C{V}_i$ such that there exist at least $r-i$ initially infected points in $ L \cap \C{V}_i$. Finally, let $E_G$ denote the event that there exist at least $r$ vertical lines passing through the $n/3 \times n/3$ grid $\C{G}$ each containing at least $r-s-1$ initially infected points in its intersection with $\C{G}$. These events are clearly independent, and it is not hard to see that the initial configuration percolates if all these events hold. Indeed, if these events all hold, then the initial configuration percolates by time $2s+3$ because at least one line (alternately from the horizontal and vertical grids) becomes active at each time $t$ with $1 \le t \le 2s+1$, and then $r$ or more vertical lines passing through $\C{G}$ become active at time $2s+2$.
Therefore,
\begin{equation}\label{lbeqn}
\vartheta_p \ge \P(E_G) \prod_{i=0}^{s}\P(E_h(i))\prod_{i=1}^{s}\P(E_v(i)).
\end{equation}

Now, for $0 \le i \le s$, given a horizontal line $L$ passing through $\C{H}_i$, the probability that there exist $r-i$ initially infected points in $L \cap \C{H}_i$ is clearly $\Theta((np/3)^{r-i})$. Since $i \le s$, we have $n(np)^{r-i} \ll 1$, so we deduce from Proposition~\ref{smallmu} that
\[\P(E_h(i)) = \Theta\left( \frac{n(np/3)^{r-i}}{3s+3}\right) = \Theta(n(np)^{r-i})\]
as $r$ and $s$ are constants not depending on $n$. Analogously, we also have $\P(E_v(i)) = \Theta(n(np)^{r-i})$ for $1 \le i \le s$. Finally, the number of vertical lines $L$ passing through $\C{G}$ such that there exist $r-s-1$ or more initially infected points in $L \cap \C{G}$ is binomially distributed with mean $\mu = \Omega(n(np)^{r-s-1})$; since $p \gg n^{-1-1/(r-s-1)}$, it follows that $\mu \gg 1$, so $\P(E_G) = \Theta(1)$ by Proposition~\ref{smallmu}. We conclude from~\eqref{lbeqn} that
\[\vartheta_p = \Omega(n^{2s+1}(np)^{r(2s+1) - s(s+1)}).\]

We now give a matching upper bound for $\vartheta_p(n,r,2)$. Given a set $A \subset [n]^2$ of initially infected points, it will be convenient to work with a \emph{slowed-down} line percolation process defined by
\[A = B^{\left(0\right)} \subset  B^{\left(1\right)} \subset \dots B^{(t)} \subset \dots \]
where
\[B^{\left(2t+1\right)} = B^{\left(2t\right)} \cup \left\{ v \in [n]^2 : \exists \text{ a horizontal line } L \in \C{L}\left(v\right) \mbox{ with } |L \cap B^{\left(2t\right)}| \ge r\right\}\]
and
\[B^{\left(2t\right)} = B^{\left(2t-1\right)} \cup \left\{ v \in [n]^2 : \exists \text{ a vertical line } L \in \C{L}\left(v\right) \mbox{ with } |L \cap B^{\left(2t-1\right)}| \ge r\right\}.\]
Since $B^{(t)} \subset A^{(t)} \subset B^{(2t)}$, it is clear that $A$ percolates if and only if $B^{(t)} = [n]^2$ for some $t\ge 0$. In this slowed-down process, a line $L \in \C{L}(n,2)$ is said to be \emph{active} at time $t$ if $L \subset B^{(t)}$ and \emph{inactive} otherwise. We say that a line $L$  \emph{fires} at time $t\ge 1$ if $L$ is active at time $t$ but not so at time $t-1$; also, a line $L$ fires at time $0$ if it is active at time $0$, i.e., if every point on $L$ is initially infected.

Let us now estimate the probability that a randomly chosen initial set $A_p \subset [n]^2$ percolates in the slowed-down line percolation process. In what follows, by the `percolation process', we shall always mean the slowed-down line percolation process described above.

First, we denote by $F_t(m)$ the event that $m$ or more lines fire at time $t$; of course, by the definition of the percolation process, these lines must all be parallel if $t \ge 1$. Next, writing $[k]_o$ and $[k]_e$ respectively for the odd and even elements of $[k]$, we say that a finite sequence $M = (m_t)_{t=1}^{k}$ of at most $2s$ positive integers is a \emph{line-count} if
\begin{enumerate}
\item $h(M) = \sum_{t \in [k]_o} m_{t} \le s$, and
\item $v(M) = \sum_{t \in [k]_e} m_t \le s$.
\end{enumerate}
For a line-count $M = (m_t)_{t=1}^{k}$, we denote by $E(M)$ the event that no lines fire at time $0$ and that exactly $m_t$ lines fire at time $t$ for all $1 \le t \le k$. 

Given a line-count $M = (m_t)_{t=1}^{k}$, we write $l(M)$ to denote $h(M)$ if $k$ is odd, and $v(M)$ if $k$ is even; clearly, $l(M) \le s$ for any line-count $M$. The following lemma is the main ingredient in the proof of the upper bound.
\begin{lemma}\label{indep}
For a line-count $M = (m_t)_{t=1}^{k}$ and any fixed integer $m \ge 0$, we have
\[ \P \left(F_{k+1}(m) \Cond E(M) \right) = O \left(\left(n\left(np\right)^{r-l(M)}\right)^{m}\right).
\]
\end{lemma}

\begin{proof}
We shall prove the lemma assuming $k$ is odd in which case $l(M) = h(M)$; the other case is analogous.

Fix a set $\C{H} \subset \C{L}(n,2)$ of $h(M)$ horizontal lines and a set $\C{V} \subset \C{L}(n,2)$ of $ v(M)$ vertical lines. Also, fix partitions $\C{H} = \C{H}_1 \cup \C{H}_3 \cup \dots \cup \C{H}_k$ and $\C{V} = \C{V}_2 \cup \C{V}_4 \cup \dots \cup \C{V}_{k-1}$ of $\C{H}$ and $\C{V}$ such that $|\C{H}_i| = m_i$ for all $i \in [k]_o$ and $|\C{V}_i| = m_i$ for all $i \in [k]_e$. Let $\Scr{E}$ denote the event that the set of lines that fire at time $t$ is precisely $\C{H}_t$ if $t \in [k]_o$ and precisely $\C{V}_t$ if $t \in [k]_e$.

Writing $E_0$ for the event that there exist no active lines at time $0$, it is clear that $E(M)$ is the disjoint union, over events $\Scr{E}$ as above, of the events $\Scr{E} \cap E_0$. Therefore, to prove the claim, it suffices to show that for each event $\Scr{E}$ as above, we have $\P (F_{k+1}(m) \Cond \Scr{E} \cap E_0) = O((n(np)^{r-l})^{m})$, where $l = l(M) = h(M)$.

Let $\C{G}$ denote the subgrid of $[n]^2$ consisting of those points not on any of the lines in  $\C{H} \cup \C{V}$. Denote by $\Scr{F}$ the event that there exist $m$ or more vertical lines passing through $\C{G}$ each containing $r-l$ or more initially infected points in its intersection with $\C{G}$. It is clear that $\P (F_{k+1}(m) \Cond \Scr{E} \cap E_0) = \P(\Scr{F}\Cond \Scr{E} \cap E_0)$. 

Now, $\Scr{F}$ is an increasing event depending only on the points in $\C{G}$ and by Proposition~\ref{smallmu}, we have $\P(\Scr{F}) = O((n(np)^{r-l})^{m})$. Also, it is clear that $E_0$ is a decreasing event. If $\Scr{E}$ was also a decreasing event, we could conclude immediately from the FKG inequality that $\P(\Scr{F}\Cond \Scr{E} \cap E_0) \le \P(\Scr{F})$. While $\Scr{E} \cap E_0$ is in itself not decreasing, we make the following useful observation.

\begin{claim} 
$\Scr{E} \cap E_0$ is decreasing on $\C{G}$.
\end{claim}
\begin{proof}
Let $\omega \in \{0,1\}^{[n]^2}$ be an initial configuration that belongs to $\Scr{E} \cap E_0$ and let $\omegap \in \{0,1\}^{[n]^2}$ be such that $\omegap_x \le \omega_x$ for $x \in \C{G}$ and $\omegap_x = \omega_x$ for $x \notin \C{G}$. We need to show that $\omegap \in \Scr{E} \cap E_0$. We first note that since the percolation process is monotone and $\omegap \le \omega$, the set of active lines at any time $t \ge 0$ when the percolation process is started from $\omegap$ is a subset of the set of active lines at that time when the percolation process is started from $\omega$. Next, note that since $\omega \in  \Scr{E} \cap E_0$, when we start the percolation process from $\omega$, none of the lines through any of the points in $\C{G}$ fire before time $k+1$; in other words, none of the points in $\C{G}$ participate in the spread of infection before time $k+1$. Therefore, since $\omegap_x = \omega_x$ for all $x \notin \C{G}$, it follows that the set of active lines at any time $t \le k$ when the percolation process is started from $\omegap$ is actually identical to the set of active lines at that time when the percolation process is started from $\omega$. It is now clear that $\omegap \in \Scr{E} \cap E_0$, as required.
\end{proof}

We now apply Proposition~\ref{fkg} to conclude that $ \P(\Scr{F}\Cond\Scr{E}\cap E_0)\le \P(\Scr{F})$, proving the lemma.
\end{proof}

Recall that $s \in\{0,1,\dots, r-1\}$. Therefore, if $A_p$ percolates, then there exists a first time $T \ge 0$ at which there exist $s+1$ or more parallel active lines; let $\HE$ denote the event that such a time $T$ exists. To bound the probability of percolation, it suffices to bound $\P (\HE)$. We do this as follows.

Let $\AM$ denote the collection of all line-counts, and for each line-count $M = (m_t)_{t=1}^{k}$ in $\AM$, denote the event $E(M) \cap F_{k+1} (s+1 -\hat l(M))$ by $\HE(M)$, where $\hat l (M) = h(M) + v(M) - l(M)$; in other words, the event $\HE(M)$ is the intersection of the event $E(M)$ and the event $\{T = k+1\}$. Also, let $\HE_0$ denote the event that there exists at least one active line at time $0$. 

Note that if $T$ exists, then $T \le 2s+1$. It follows that
\[\HE \subset \bigcup_{M \in \AM}\HE(M) \cup \HE_0\]
and consequently,
\[\vartheta_p \le \P(\HE) \le \sum_{M \in \AM}\P(\HE(M)) + \P(\HE_0).\]

First, we bound $\P(\HE(M))$ for each $M \in \AM$. By Lemma~\ref{indep}, for a line-count $M = (m_t)_{t=1}^{k}$ in $\AM$, we have
\begin{align*}
P(\HE(M)) &= \P(E(M))\P\left(F_{k+1} (s+1 - \hat l (M)) \Cond E(M)\right) \\
&= O\left( \P(E(M))\left(n(np)^{r-l(M)}\right)^{s+1 - \hat l (M)}\right).
\end{align*}
If $M$ is the empty sequence, then we trivially have $\P(E(M)) \le 1$. If $M = (m_t)_{t=1}^{k}$ with $k \ge 1$ on the other hand, then let $M' = (m_t)_{t=1}^{k-1}$ and note that by Lemma~\ref{indep}, we have
\begin{align}
\P(E(M)) &= \P(E(M')) \P\left(F_k(m_k) \setminus F_k (m_k + 1) \Cond E(M')\right) \nonumber  \\
&\le \P(E(M')) \P(F_k(m_k)\Cond E(M')) \nonumber \\
&= O\left( \P(E(M')) \left(n(np)^{r-l(M')}\right)^{m_k}\right).\label{ind-linecount}
\end{align}
Since the length of any line-count is at most $2s$ (which is a constant not depending on $n$), we may bound $P(\HE(M))$ inductively using~\eqref{ind-linecount}; indeed, it follows from $k \le 2s$ applications of~\eqref{ind-linecount} that
\begin{align*}
\P(\HE(M)) = O\Bigl( &\left(n(np)^r\right)^{m_1} \times \left(n(np)^{r-m_1}\right)^{m_2} \times \left(n(np)^{r-m_2}\right)^{m_3} \\
& \times \left(n(np)^{r-m_1-m_3}\right)^{m_4} \times \dots \times \left(n(np)^{r-l(M)}\right)^{s+1 - \hat l(M)} \Bigr).
\end{align*}
This bound, on algebraic simplification, shows that for any $M \in \AM$, we have
\begin{align*}
\P(\HE(M))&= O\left(n^{s+1+l(M)}(np)^{r(s+1) + (r-s-1)l(M)}\right)\\
&= O\left(n^{s+1}(np)^{r(s+1)} \left(n(np)^{r-s-1}\right)^{l(M)} \right)\\
&= O\left(n^{2s+1}(np)^{r(s+1) + (r-s-1)s}\right)\\
&= O\left(n^{2s+1}\left(np\right)^{r(2s+1) - s(s+1)}\right);
\end{align*}
here, we have used the fact that $n(np)^{r-s-1} \gg 1$ and $l(M) \le s$.

Next, by a simple union bound, we have $\P(\HE_0) \le 2np^n$. Furthermore, it is easy to check that $2np^n = o(n^{2s+1}\left(np\right)^{r(2s+1) - s(s+1)})$ for any $p \le 1/n$.

Finally, we know from the definition of a line-count that each element of a line-count is a positive integer not exceeding $s$, and that the length of a line-count is at most $2s$. Consequently, the number of line-counts is at most $s^{2s}$. It is now clear that
\begin{align*}
\vartheta_p \le \P(\HE) &=  O\left(|\AM | n^{2s+1}\left(np\right)^{r(2s+1) - s(s+1)}\right)\\
&= O\left( n^{2s+1}\left(np\right)^{r(2s+1) - s(s+1)}\right). \qedhere
\end{align*}
\end{proof}

\section{The critical probability in three dimensions}\label{3d}

We now turn our attention to the line percolation process in three dimensions and prove Theorem~\ref{3d-perc}.

\begin{proof}[Proof of Theorem~\ref{3d-perc}]

Suppose that the points of $[n]^3$ are initially infected independently with probability $p$ and let $C = C(n) = p/n^{-1-1/(r-\gamma)}$. We shall show that percolation occurs with probability at least $1/2$ provided $C$ is greater than some sufficiently large constant, and that percolation occurs with probability at most $1/2$ provided $C$ is less than some sufficiently small constant. For technical reasons, we shall establish this under the additional assumption that $C$ is neither too large nor too small; of course, we are free to do this in the light of monotonicity, so in what follows, we shall assume that $1/\log n \le C \le \log n$ without any loss of generality. We prove the upper and lower bounds separately by distinguishing the following two cases.

\textbf{Case 1: $C \gg 1$.} Unsurprisingly, it is easier to show that percolation occurs than to demonstrate otherwise. We start by bounding $p_c$ from above by showing that percolation occurs with probability at least $1/2$ provided $C$ is greater than some sufficiently large constant. Note that $s = \lfloor \sqrt{r+ 1/4} - 1/2\rfloor$, as defined in the statement of Theorem~\ref{3d-perc}, is the unique positive integer such that 
\[s(s+1) \le r < (s+1)(s+2).\]
It is not hard to check from the definition of $s$ that $\gamma = (r + s^2 + s)/(2s+2)$ satisfies
\[
n^{-1-\frac{1}{r-s-1}} \ll n^{-1-\frac{1}{r-\gamma}}\ll n^{-1-\frac{1}{r-s}},
\]
so it follows from Theorem~\ref{2d-unboosted-perc} that
\[\vartheta_p\left(n,r,2\right) = \Theta\left(n^{2s+1}\left(np\right)^{r\left(2s+1\right)-s\left(s+1\right)}\right) = \Theta\left(C^{r\left(2s+1\right)-s\left(s+1\right)}n^{-1}\right),\]
where the implicit constants suppressed by the asymptotic notation depend only on $r$ and $s$.

We say that a plane $P$ is \emph{internally spanned} if $A^{(0)}\cap P$ percolates in the line percolation process restricted to $P$. Choose an axis-parallel direction and consider the $n$ parallel planes perpendicular to this direction. The number of such planes  that are internally spanned is a binomial random variable with mean $\mu = \Omega(C^{r(2s+1)-s(s+1)})$. Since $\mu \rightarrow \infty$ as $C \rightarrow \infty$, we see from Proposition~\ref{chernoff} that there exist $r$ parallel internally spanned planes with probability at least $1/2$ if $C$ exceeds a sufficiently large constant. It follows that $p_c(n,r,3) = O(n^{-1-1/(r-\gamma)})$.

\textbf{Case 2: $C \ll 1$.} We claim that the probability of percolation is at most $1/2$ provided $C$ is less than some sufficiently small constant. It will be helpful to handle the case where $r = 2$ separately. Indeed, if $r = 2$ (and $s = \gamma = 1$), then the expected number of lines $L \in \C{L}(n,3)$ that contain $r$ initially infected points is $O(n^2 (np)^r) = O(n^4 p^2) = O(C^2)$, and consequently, there are no such lines with probability at least $1/2$ if $C$ is less than some suitably small constant, implying that $p_c(n,2,3) = \Omega(n^{-2})$. In the sequel, we therefore suppose that $r \ge 3$.

We shall demonstrate that the probability of percolation is at most $1/2$ by proving something much stronger. We shall track, as the infection spreads, the number of planes containing $\ell$ or more parallel active lines for each $1 \le \ell \le s+1$ and show that, with probability at least $1/2$, these numbers are not too large when the process terminates; in particular, we shall show that there are no planes with $s+1$ or more parallel active lines when the process terminates and consequently, that there is no percolation.

As we did in the two-dimensional case, we shall work with a modified three-dimensional line percolation process in which the infection only spreads along a single direction at each stage. For a set $A \subset [n]^3$ of initially infected points, denoting the standard basis for $\R^3$ by $\{e_1, e_2, e_3\}$, in the \emph{modified} line percolation process, we have a sequence
\[A = B^{(0)} \subset B^{(1)} \subset \dots \subset B^{(t)} \subset \dots\] of subsets of $[n]^3$ where $B^{(t+1)}$ is obtained from $B^{(t)}$ by spreading the infection only along lines parallel to $e_i$ where $i \equiv t\imod{3}$. Furthermore, we terminate this modified process at the first time $t\ge 0$ at which
\begin{enumerate}
\item either $B^{(t)} = [A]$, or
\item \label{condA} for some $1 \le \ell \le s+1$, there exist $n^{1-\ell\gamma/(r-\gamma)}$ or more planes that each contain $\ell$ or more parallel active lines.
\end{enumerate}
We say that $A$ \emph{diverges} if condition~(\ref{condA}) holds when the modified process terminates. It is clear that if $A$ percolates, then $A$ diverges in the modified process. 

We need some definitions analogous to those in the proof of Theorem~\ref{2d-unboosted-perc}. In the modified process, a line $L \in \C{L}(n,3)$ is said to be \emph{active} at time $t$ if $L \subset B^{(t)}$ and \emph{inactive} otherwise. As before, a line $L$  \emph{fires} at time $t\ge 1$ if $L$ is active at time $t$ but not so at time $t-1$; also, a line $L$ fires at time $0$ if it is active at time $0$. 

In the sequel, by the `percolation process', we shall always mean the modified line percolation process described above. We shall estimate the probability that a randomly chosen initial set $A_p \subset [n]^3$ diverges in the percolation process.  Writing $D$ to denote the event that $A_p$ diverges in the percolation process, we prove the following bound.

\begin{claim}\label{planecount} $\P\left(D\right) = O\left(\sum_{\ell = 1}^s C^{r\ell} + C^{r\left(2s+1\right)-s\left(s+1\right)}\right).$
\end{claim}

Before we prove Claim~\ref{planecount}, let us note that the required lower bound for the critical probability follows immediately from the claim. Indeed, the claim implies that $\P(D) \rightarrow 0$ as $C \rightarrow 0$ for all $r \ge 3$, implying that $p_c(n,r,3) = \Omega(n^{-1-1/(r-\gamma)})$.

\begin{proof}[Proof of Claim~\ref{planecount}]
Fix a plane $P$ and let $\C{L}(P)$ denote the set of $2n$ lines contained in $P$. We shall prove Claim~\ref{planecount} by estimating the probability that $P$ contains $\ell$ or more parallel active lines when the percolation process terminates. The spread of infection within $P$ resembles the two-dimensional line percolation process, with the key distinction that some points in $P$ also become infected by virtue of lying on an active line perpendicular to $P$. However, since we are interested in estimating the probability that $\ell$ or more parallel lines in $P$ fire \emph{before $A_p$ diverges}, we shall not have to worry about there being too many such points.

For $1 \le \ell \le s+1$, let $Q_\ell$ denote the event that $\ell$ or more parallel lines in $\C{L}(P)$ are active when the percolation process terminates. We prove the following claim.

\begin{claim}\label{PE_k}
For $1 \le \ell \le s+1$, we have
\[\P\left(Q_\ell\right) = 
\begin{cases}
O\left( C^{r\ell} n^{-\ell\gamma / \left(r-\gamma\right)}\right) \hfill & \text{ if } 1 \le \ell \le s, \text{ and}\\
O\left( C^{r\left(2s+1\right)-s\left(s+1\right)} n^{-1}\right) \hfill & \text{ if } \ell = s+1.
\end{cases}
 \]
\end{claim}
\begin{proof}
Assume for concreteness that $P$ is perpendicular to $e_3$ in which case each line in $\C{L}(P)$ is parallel to either $e_1$ or $e_2$; we shall think of the lines in $\C{L}(P)$ parallel to $e_1$ as being `horizontal' and the lines parallel to $e_2$ as being `vertical'.

Let us fix $\ell \in [s+1]$. As in the two-dimensional case, we shall bound $\P(Q_\ell)$ by estimating the probability of this event happening according to a particular `line-count'. Note that unlike in the two-dimensional process, a large amount of time may elapse between two successive stages at which lines in $\C{L}(P)$ fire. Consequently, the precise notion of a line-count that we use here differs slightly from the notion used previously.

We call a time $t \ge 0$ an \emph{epoch (for $P$)} if at least one line in $\C{L}(P)$ fires at time $t$. We denote by $H_i(m)$ the event that $m$ or more horizontal lines in $\C{L}(P)$ fire in the $i$th epoch, and we define the event $V_i(m)$ analogously.

A \emph{line-count} $M = ((m_i, d_i))_{i = 1}^k$ is a sequence of at most $2(\ell - 1)$ pairs $(m_i, d_i)$ with $m_i\in\N$ and $d_i\in \{e_1,e_2\}$ such that
\begin{enumerate}
\item $h(M) = \sum_{i:d_i = e_1} m_i < \ell$, and
\item $v(M) = \sum_{i:d_i = e_2} m_i < \ell$.
\end{enumerate}
For a line-count $M = ((m_i, d_i))_{i=1}^k$, let $E(M)$ denote the event that no lines in $\C{L}(P)$ are active at time $0$ and that the number and direction of the lines that fire in the $i$th epoch are given by $m_i$ and $d_i$ respectively for all $1 \le i \le k$.

Note that since $\ell \in [s+1]$, both $h(M)$ and $v(M)$ are at most $s$ for any line-count $M$. The key ingredient in the proof of the upper bound is the following lemma.

\begin{lemma}\label{indep'}
For any line-count $M = ((m_i, d_i))_{i = 1}^k$ and any fixed integer $m \ge 0$, we have
\[ \P \left(V_{k+1}(m) \Cond E(M) \right) = O \left(\left(n\left(np\right)^{r-h(M)}\right)^{m}\right)
\]
and similarly,
\[ \P \left(H_{k+1}(m) \Cond E(M) \right) = O \left(\left(n\left(np\right)^{r-v(M)}\right)^{m}\right).
\]
\end{lemma}

\begin{proof}
We only prove the first claim; the proof of the other assertion is analogous. The argument we adopt is very similar to the one used to prove Lemma~\ref{indep}. However, the main difference is that we need to account for points of $P$ that become infected between two epochs by virtue of lying on a line perpendicular to $P$ that fires between epochs. We call a point of $P$ a \emph{boost} if the line perpendicular to $P$ through that point fires before the $(k+1)$th epoch.

Fix disjoint sets of lines $\C{L}_1, \C{L}_2, \dots, \C{L}_k \subset \C{L}(P)$ so that the set $\C{L}_i$ consists of $m_i$ lines parallel to $d_i$ for each $1 \le i \le k$. Let $\Scr{E}'$ denote the event that the set of lines infected in the $i$th epoch is precisely $\C{L}_i$ for all $1 \le i \le k$. Also, fix a set $\C{B} \subset P$ of points and let $\Scr{E}''$ denote the event that set of boosts in $P$ is precisely $\C{B}$. 

Writing $E_0$ for the event that there exist no active lines in $\C{L}(P)$ at time $0$, it is clear that $E(M)$ is the disjoint union, over events $\Scr{E}'$ as above, of the events $\Scr{E}' \cap E_0$. Therefore, to prove the claim, it suffices to show that for each event $\Scr{E}'$ as above, we have $\P (V_{k+1}(m) \Cond \Scr{E}' \cap E_0) = O((n(np)^{r-h(M)})^{m})$. We shall demonstrate the stronger assertion that 
\[\P (V_{k+1}(m) \Cond \Scr{E}' \cap \Scr{E}'' \cap E_0) = O\left(\left(n(np)^{r-h(M)}\right)^{m}\right)\]
for every pair of events $\Scr{E}'$ and $\Scr{E}''$ as above. 

For simplicity of notation, we denote the event $\Scr{E}' \cap \Scr{E}'' \cap E_0$ by $\Scr{E}$ and write $h = h(M)$ and $v = v(M)$.

Denote by $\C{V}_j$ the set of those vertical lines of $\C{L}(P)$ that meet $\C{B}$ in $j$ points and let $N_j = |\C{V}_j|$. When we condition on $\Scr{E}''$, note that we assume that there exist at least $N_j$ planes containing $j$ or more parallel active lines before the $(k+1)$th epoch. Therefore, if $N_j > n^{1-j\gamma/(r-\gamma)}$ for some $1 \le j \le s+1$, then 
\[\P (V_{k+1}(m) \Cond \Scr{E}) = 0;\] 
indeed, $V_{k+1}(m)$ and $\Scr{E}''$ are disjoint in this case since conditioning on $\Scr{E}''$ tells us that the percolation process terminates before a $(k+1)$th epoch may even occur. Hence, we may assume that $N_j \le n^{1-j\gamma/(r-\gamma)}$ for each $1 \le j \le s+1$. Observe that $(s+1)\gamma / (r-\gamma) > 1$ since $(s+1)(s+2) > r$, so $N_j = 0$ for $j \ge s+1$ as $n^{1-(s+1)\gamma/(r-\gamma)} < 1$.

Let $\C{G}$ denote the subgrid of $[n]^2$ consisting of those points not on any of the lines in $\C{L}_1  \cup \C{L}_2 \cup \dots \cup \C{L}_k$. For $0 \le j \le s$, we write $X_j$ for the (random) number of lines in $\C{V}_j$ whose intersection with $\C{G} \setminus \C{B}$ contains $r-h-j$ or more initially infected points. Let $\Scr{F}$ denote the event that $\sum_{j= 0}^{s} X_j \ge m$. Note that a point of $\C{G}$ that is infected before the $(k+1)$th epoch is either initially infected or belongs to $\C{B}$; therefore, 
\[\P (V_{k+1}(m) \Cond \Scr{E}) \le \P ( \Scr{F} \Cond \Scr{E}). \]

\begin{figure}
\begin{center}
\begin{tikzpicture}[xscale=1.6,yscale=1]
\draw (0,0)--(5,0)--(7,2.5)--(2,2.5)--(0,0);
\draw (1.3,4.5)--(6.3,4.5)--(6.3,-0.5)--(1.3,-0.5)--(1.3,4.5);
\draw [dotted, thick] (1.3,1.625)--(6.3,1.625);
\draw [very thick](2.5,-0.5)--(2.5,4.5);
\draw [very thick](3.4,-0.5)--(3.4,4.5);
\draw [very thick](5.4,-0.5)--(5.4,4.5);
\node (plane) at (0.6, 1.3) {$\C{G}$};
\node (b1) at (2.5, 1.625) [inner sep=0.9mm, thick, rectangle, draw=black!100, fill=black!50] {};
\node (b2) at (3.4, 1.625) [inner sep=0.9mm, thick, rectangle, draw=black!100, fill=black!50] {};
\node (b3) at (5.4, 1.625) [inner sep=0.9mm, thick, rectangle, draw=black!100, fill=black!50] {};
\node (b1-1) at (3.4, 3.825) [inner sep=0.7mm, thick, circle, draw=black!100, fill=black!50] {};
\node (b1-2) at (3.4, 3.475) [inner sep=0.7mm, thick, circle, draw=black!100, fill=black!50] {};
\node (b1-3) at (3.4, 2.875) [inner sep=0.7mm, thick, circle, draw=black!100, fill=black!50] {};
\node (b1-4) at (3.4, 2.25) [inner sep=0.7mm, thick, circle, draw=black!100, fill=black!50] {};
\node (b2-1) at (2.5, 3.9) [inner sep=0.7mm, thick, circle, draw=black!100, fill=black!50] {};
\node (b2-2) at (2.5, 3.5) [inner sep=0.7mm, thick, circle, draw=black!100, fill=black!50] {};
\node (b2-3) at (2.5, 3.1) [inner sep=0.7mm, thick, circle, draw=black!100, fill=black!50] {};
\node (b2-4) at (2.5, 0.3) [inner sep=0.7mm, thick, circle, draw=black!100, fill=black!50] {};
\node (b3-1) at (5.4, 2.7) [inner sep=0.7mm, thick, circle, draw=black!100, fill=black!50] {};
\node (b3-2) at (5.4, 2.3) [inner sep=0.7mm, thick, circle, draw=black!100, fill=black!50] {};
\node (b3-3) at (5.4, 0.9) [inner sep=0.7mm, thick, circle, draw=black!100, fill=black!50] {};
\node (b3-1) at (5.4, -0.1) [inner sep=0.7mm, thick, circle, draw=black!100, fill=black!50] {};
\end{tikzpicture}
\end{center}
\caption{Boosted-points on a line in $\C{G}$.}
\label{boostpic}
\end{figure}
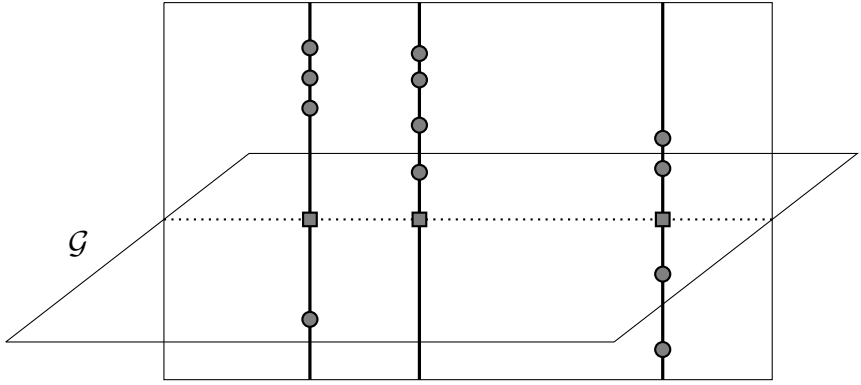

Clearly, $\Scr{F}$ is an increasing event that depends only on the points in $\C{G} \setminus \C{B}$. As before, we have the following observation.

\begin{claim}
$\Scr{E}$ is decreasing on $\C{G} \setminus \C{B}$. 
\end{claim}
\begin{proof}
Let $\omega \in \{0,1\}^{[n]^3}$ be an initial configuration that belongs to $\Scr{E}$ and let $\omegap \in \{0,1\}^{[n]^3}$ be such that $\omegap_x \le \omega_x$ for $x \in \C{G} \setminus \C{B}$ and $\omegap_x = \omega_x$ for $x \notin \C{G} \setminus \C{B}$. We need to show that $\omegap \in \Scr{E}$. As before, we first note that since the percolation process is monotone and $\omegap \le \omega$, the set of active lines (not just in $P$, but rather in all of $[n]^3$) at any time $t \ge 0$ when we start from $\omegap$ is a subset of the set of active lines at that time when we start from $\omega$. Next, note that since $\omega \in \Scr{E}$, when the percolation process is started from $\omega$, none of the lines through any of the points in $\C{G}\setminus \C{B}$ fire before the $(k+1)$th epoch; in other words, none of these points participate in the spread of infection before the $(k+1)$th epoch. Therefore, since $\omegap_x = \omega_x$ for all $x \notin \C{G} \setminus \C{B}$, it follows that the set of active lines at any time before the $(k+1)$th epoch when we start from $\omegap$ is actually identical to the set of active lines at that time when we start from $\omega$. It follows that $\omegap \in \Scr{E}$, implying that $\Scr{E}$ is decreasing on $\C{G} \setminus \C{B}$.
\end{proof}

It now follows from Proposition~\ref{fkg} that $\P(\Scr{F} \Cond \Scr{E}) \le \P(\Scr{F})$. Therefore, to finish the proof of the lemma, it suffices to show that $\P(\Scr{F}) = O((n(np)^{r-h})^{m})$.

For $0 \le j \le s$, consider the event $\{X_j \ge x\}$, i.e., the event that there exist $x$ or more vertical lines in $\C{V}_j$ that each meet $\C{G} \setminus \C{B}$ in $r-h-j$ or more initially infected points. It is easy to see from Proposition~\ref{smallmu} that
\[ \P(X_0 \ge x ) = O \left( \binom{N_0}{x}\left(\left(np\right)^{r-h}\right)^{x}\right) = O\left(\left(n\left(np\right)^{r-h}\right)^{x}\right),\]
where the last inequality holds since $N_0 \le n$. For $1 \le j \le s$ on the other hand, we have
\begin{align*}
 \P (X_j \ge x) &= O \left( \binom{N_j}{x}\left(\left(np\right)^{r-h-j}\right)^{x}\right)\\
 &= O\left(\left(n\left(np\right)^{r-h}\right)^{x} \left(n^{1 + \gamma/\left(r-\gamma\right)}p\right)^{-jx}\right)\\
 &= O\left(\left(n\left(np\right)^{r-h}\right)^{x}\right)
\end{align*}
since $N_j \le n^{1-j\gamma/(r-\gamma)}$ and $n^{1+\gamma/(r-\gamma)}p=Cn^{(\gamma - 1)/(r-\gamma)}$, noting that $ Cn^{(\gamma - 1)/(r-\gamma)} = \omega(1)$ since $C \ge 1/\log n$ and $\gamma > 1$ for all $r \ge 3$. It follows that
\[ \P(\Scr{F}) = \P\left(\sum_{j=0}^s X_j \ge m\right) \le \sum_{x_0, x_1, \dots,x_s} \P(X_0 \ge x_0, X_1 \ge x_1, \dots, X_s \ge x_s),
\]
where the sum above is over all non-negative integer solutions to the equation $x_0 + x_1 + \dots +x_s = m$. First, note that the number of such solutions is at most $(m+1)^s$ which is a constant not depending on $n$ since $m$ is fixed. Next, note that since the sets $\C{V}_0, \C{V}_1, \dots, \C{V}_s$ are disjoint, the random variables $X_0, X_1, \dots, X_s$ are independent. Therefore,
\[\P(X_0 \ge x_0, X_1 \ge x_1, \dots, X_s \ge x_s) = \prod_{j=0}^s \P(X_j \ge x_j) = O\left(\left(n\left(np\right)^{r-h}\right)^{\sum_{j=0}^s x_j}\right).
\]
It follows that $\P(\Scr{F}) \le O((n(np)^{r-h})^{m})$, proving the lemma.
\end{proof}

Let $\AM$ denote the collection of all line-counts. For a line-count $M = ((m_i,d_i))_{i=1}^{k}$ in $\AM$, let 
\[\HE(M) = E(M) \cap (V_{k+1} (\ell - v(M)) \cup H_{k+1} (\ell - h(M))).\] Also, let $\HE_0$ denote the event that there exists at least one active line at time $0$. As in the proof of Theorem~\ref{2d-unboosted-perc}, we have
\[Q_\ell \subset \bigcup_{M \in \AM}\HE(M) \cup \HE_0\]
and consequently,
\[\P(Q_\ell) \le \sum_{M \in \AM}\P(\HE(M)) + \P(\HE_0).\]

First, we bound $\P(\HE(M))$ for each $M \in \AM$. Fix a line-count $M = ((m_i,d_i))_{i=1}^{k}$ in $\AM$ and let $h = h(M)$ and $v = v(M)$. Clearly, we have
\[
P(\HE(M)) \le \P(E(M))\left(\P\left(V_{k+1} (\ell - v) \Cond E(M)\right) + \P\left(H_{k+1} (\ell - h) \Cond E(M)\right) \right).
\]
As $\ell \le s+1$ is a constant not depending on $n$, by $k+1 \le 2\ell - 1$ applications of Lemma~\ref{indep'}, we get
\[
\P(\HE(M)) = O\left(\max\left((n(np)^{r-h})^{\ell-v}, (n(np)^{r-v})^{\ell-h}\right) \prod_{i=1}^{k}\left(n(np)^{r-\sum_{j < i: d_j \neq d_i}m_j}\right)^{m_{i}}\right).
\]
The above bound, on algebraic simplification, yields
\[
\P(\HE(M)) = O\left(\max\left( n^{\ell+h}\left(np\right)^{r\ell + rh - \ell h} , n^{\ell+v}\left(np\right)^{r\ell + rv - \ell v}\right)\right).
\]
Since $0 \le h(M), v(M) < \ell$ for any $M \in \AM$, it follows that
\[
\P(\HE(M)) = O \left( \max_{0 \le x < \ell} \left( n^{\ell+x}\left(np\right)^{r\ell + rx - \ell x} \right) \right) = O\left(\max_{0 \le x < \ell}\left(n^{\ell}\left(np\right)^{r\ell}\left(n\left(np\right)^{r-\ell}\right)^{x}\right)\right)
\]
for any $M \in \AM$. 

Next, by a simple union bound, we have $\P(\HE_0) \le 3n^2p^n$. Furthermore, it is easy to check that $3n^2p^n = o(n^\ell (np)^{r\ell})$ for any $p \le 1/n$.

Since the number of line-counts is clearly at most  $\ell^{2\ell}$ (which is a constant not depending on $n$), it is now easy to check that
\begin{equation}
\P(Q_\ell) = O\left(\max_{0 \le x < \ell} \left(n^{\ell}\left(np\right)^{r\ell}\left(n\left(np\right)^{r-\ell}\right)^{x} \right)\right). \label{Ek-prob}
\end{equation}
Recall that $p = Cn^{-1-1/(r-\gamma)}$. It is easy to check from the definitions of $s$ and $\gamma$ that the estimate for the probability of $Q_\ell$ in~\eqref{Ek-prob} is maximised by taking $x = 0$ when $1 \le \ell \le s$, and by taking $x = s$ when $\ell = s+1$. Consequently,
\[\P\left(Q_\ell\right) = 
\begin{cases}
O\left(n^\ell (np)^{r\ell}\right) \hfill & \text{ if } 1 \le \ell \le s, \text{ and}\\
O\left(n^{2s+1}\left(np\right)^{r\left(2s+1\right)-s\left(s+1\right)}\right) \hfill & \text{ if } \ell = s+1,
\end{cases}
\]
which on algebraic simplification yields
\[
\P\left(Q_\ell\right) = 
\begin{cases}
O\left( C^{r\ell} n^{-\ell\gamma / \left(r-\gamma\right)}\right) \hfill & \text{ if } 1 \le \ell \le s, \text{ and}\\
O\left( C^{r\left(2s+1\right)-s\left(s+1\right)} n^{-1}\right) \hfill & \text{ if } \ell = s+1,
\end{cases}
\]
proving the claim.
\end{proof}

Recall that $D$ is the event that when the percolation process terminates, the number of planes containing $\ell$ or more parallel active lines exceeds $  n^{1-\ell\gamma/(r-\gamma)} $ for some $1\le \ell \le s+1$.

From Claim~\ref{PE_k}, we see that the expected number of planes containing at least $\ell$ parallel active lines when the percolation process terminates is $O( C^{r\ell} n^{1-\ell\gamma/(r-\gamma)} )$ when $1\le \ell \le s$ and $O(C^{r(2s+1)-s(s+1)})$ when $\ell=s+1$. By Markov's inequality, the probability that the number of planes containing $\ell$ or more parallel active lines exceeds $n^{1-\ell\gamma/(r-\gamma)} $ is $O(C^{r\ell})$ when $1\le \ell \le s$ and $O(C^{r(2s+1)-s(s+1)})$ when $\ell=s+1$ since $ n^{1-(s+1)\gamma/(r-\gamma)} < 1$. Applying the union bound, we get
\[ \P\left(D\right) = O\left(\sum_{\ell=1}^s C^{r\ell} + C^{r\left(2s+1\right)-s\left(s+1\right)}\right).\qedhere\]
\end{proof}

We have now established Claim~\ref{planecount}; as described earlier, Theorem~\ref{3d-perc} is an immediate consequence.
\end{proof}	

\section{Minimal percolating sets}\label{minset}

In this section, we prove Theorem~\ref{min-set} which tells us the size of a minimal percolating set. We shall make use of the polynomial method which has had many surprising applications in combinatorics; see~\citep{Guth13} for a survey of many of these applications. While linear algebraic techniques have previously been used to study bootstrap percolation processes (see~\citep{linalg}), we believe that this application of the polynomial method is new to the field.

\begin{proof}[Proof of Theorem~\ref{min-set}]
Suppose for the sake of contradiction that there is a set $A\subset [n]^d$ which percolates with $|A| < r^d$. We shall derive a contradiction using the polynomial method.

\begin{proposition}
There exists a nonzero polynomial $P_A \in \R[x_1, x_2, \dots, x_d]$ of degree at most $r-1$ in each variable that vanishes on $A$.
\end{proposition}
\begin{proof}
Let $V \subset \R[x_1, x_2, \dots, x_d]$ be the vector space of real polynomials in $d$ variables of degree at most $r-1$ in each variable. Clearly, $\dim(V) = r^d$. Consider the evaluation map from $V$ to $\R^{|A|}$ which sends a polynomial $P$ to $(P(v))_{v \in A}$. This map is linear, and since we assumed that $|A| < r^d$, this map has a nontrivial kernel. The existence of $P_A$ follows.
\end{proof}

We shall use the polynomial $P_A$ to follow the spread of infection. The following claim will yield a contradiction.

\begin{proposition}
The polynomial $P_A$ vanishes on $A^{(t)}$ for every $t\ge 0$.
\end{proposition}

\begin{proof}
We proceed by induction on $t$. The claim is true when $t = 0$ since $A^{(0)} = A$. Now, assume $P_A$ vanishes on $A^{(t)}$ and consider a line $L$ that is active at time $t+1$ but inactive at time $t$. It must be the case that $|L \cap A^{(t)}| \ge r$. Since  $P_A$ vanishes on $A^{(t)}$, the restriction of $P_A$ to $L$ disappears on $L \cap A^{(t)}$. Denoting the standard basis for $\R^d$ by $\{e_1, e_2, \dots, e_d\}$, note that if $L$ is parallel to $e_i$, then the restriction of $P_A$ to $L$ is a univariate polynomial in the variable $x_i$ of degree at most $r-1$. Since a nonzero univariate polynomial of degree at most $r-1$ has no more than $r-1$ roots, the restriction of $P_A$ to $L$ has to be identically zero. Consequently, $P_A$ vanishes on $A^{(t+1)}$.
\end{proof}

Since $A$ percolates, we conclude that $P_A$ vanishes on $[n]^d$. On the other hand, using the following proposition, the proof of which may be found in~\citep{Alon99}, we conclude that $P_A$ cannot vanish on $[r+1]^d$.

\begin{proposition}\label{weak-null} Let $P = P(x_1, x_2, \dots, x_d)$ be a polynomial in $d$ variables over an arbitrary field $\Scr{F}$. For $1\le i \le d$, suppose that the degree of $P$ as a polynomial in $x_i$ is at most $k_i$, and let $S_i \subset \Scr{F}$ be a set of size at least $k_i + 1$. If $P(u_1, u_2, \dots, u_d) = 0$ for every $d$-tuple $(u_1, u_2, \dots, u_d) \in S_1 \times S_2 \times \dots \times S_d$, then $P$ is identically zero. \qed
\end{proposition}

It follows from Proposition~\ref{weak-null} that $P_A$ is identically zero, which is a contradiction. This establishes Theorem~\ref{min-set}.
\end{proof}

\begin{remark}
It follows from Theorem~\ref{min-set} that the size of a minimal percolating set in the $r$-neighbour bootstrap percolation model on $[n]^d$ with edges induced by the Hamming torus is at least $(r/d)^d$. On the other hand, it is possible to construct sets of size about $r^d/2d!$ which percolate. It would be interesting to determine the size of a minimal percolating set in this model exactly for all $d,r \in \N$; we suspect that the lower bound of $(r/d)^d$ is quite far from the truth.
\end{remark}

\section{Conclusion}\label{remark}
There remain many challenging and attractive open problems, chief amongst which is the determination of $p_c(n,r,d)$ for all $d, r \in \N$. To determine $p_c(n,r,3)$, we used a careful estimate for $\vartheta_p(n,r,2)$ which is valid for all $0 \le p \le 1$. This estimate for $\vartheta_p(n,r,2)$ depends crucially on the fact that the two-dimensional process reaches termination in a bounded number (depending on $r$, but not on $n$) of steps. We believe that to determine $p_c(n,r,4)$, one will need to determine $\vartheta_p(n,r,3)$ for all $0 \le p \le 1$ but since it is not at all obvious that the three-dimensional process reaches termination in a bounded number of steps with high probability, we suspect different methods will be necessary.

As remarked earlier, it is easily read out of our proofs that the line percolation model does not have a sharp threshold in two or three dimensions. It would be interesting to show that an analogous statement holds in all dimensions.

\section*{Acknowledgements}
The first and second authors were partially supported by NSF grant DMS-1301614 and the second author also wishes to acknowledge support from EU MULTIPLEX grant 317532.

Some of the research in this paper was carried out while the authors were visitors at Microsoft Research, Redmond. This research was continued while the third and fourth authors were visitors at the University of Memphis. The authors are grateful to Yuval Peres and the other members of the Theory Group at Microsoft Research for their hospitality, and the third and fourth authors are additionally grateful for the hospitality of the University of Memphis.

We would also like to thank Oliver Riordan and the anonymous referees for their helpful comments.

\bibliographystyle{amsplain}
\bibliography{line_perc}

\end{document}